\documentclass[11pt]{article}
\usepackage{amssymb,latexsym,geometry,cite,enumerate,float}
\newtheorem{theorem}{Theorem}

\newtheorem{corollary}[theorem]{Corollary}

\usepackage{setspace}
\usepackage{lineno}

\begin{document}

\onehalfspace

\title{Forbidden Induced Subgraphs for Bounded $p$-Intersection Number}

\author{Claudson F. Bornstein$^1$, Jos\'{e} W.C. Pinto$^1$, Dieter Rautenbach$^2$, Jayme L. Szwarcfiter$^1$}

\date{}

\maketitle

\begin{center}
$^1$ 
Federal University of Rio de Janeiro, Rio de Janeiro, Brazil,\\
\texttt{cbornstein@gmail.com, jwcoura@cos.ufrj.br, jayme@nce.ufrj.br}\\[3mm]
$^2$ 
Institute of Optimization and Operations Research, Ulm University, Ulm, Germany,
\texttt{dieter.rautenbach}\texttt{@uni-ulm.de}
\end{center}

\begin{abstract}
A graph $G$ has $p$-intersection number at most $d$
if it is possible to assign to every vertex $u$ of $G$, a subset $S(u)$ of some ground set $U$ with $|U|=d$
in such a way that distinct vertices $u$ and $v$ of $G$ are adjacent in $G$ if and only if $|S(u)\cap S(v)|\geq p$.
We show that every minimal forbidden induced subgraph 
for the hereditary class ${\cal G}(d,p)$ 
of graphs whose $p$-intersection number is at most $d$, 
has order at most $3\cdot 2^{d+1}+1$,
and that the exponential dependence on $d$ in this upper bound is necessary.
For $p\in \{ d-1,d-2\}$, we provide more explicit results 
characterizing the graphs in ${\cal G}(d,p)$ 
without isolated/universal vertices
using forbidden induced subgraphs.
\end{abstract}


\section{Introduction}

Intersection representations of graphs are among the most important graph representations
and lead to some famous and well studied graph classes \cite{mm}.
As a generalization of intersection representations, 
Jacobson et al. \cite{jms} introduced $p$-intersection representations.
For a positive integer $p$, a {\it $p$-intersection representation} of a graph $G$ is a function $S:V(G)\to 2^U$ 
assigning to every vertex $u$ of $G$, a subset $S(u)$ of some ground set $U$ 
in such a way that 
distinct vertices $u$ and $v$ are adjacent in $G$ 
if and only if $|S(u)\cap S(v)|\geq p$.
The choice $p=1$ leads to classical intersection representations of graphs.
Since every graph has a $p$-intersection representation for every $p$,
it makes sense to study the {\it $p$-intersection number $\Theta_p(G)$} of $G$
defined as the minimum cardinality of a set $U$ for which $G$ has a 
$p$-intersection representation $S:V(G)\to 2^U$ with ground set $U$.
The $1$-intersection number was first studied by Erd\H{o}s et al. \cite{egp}
who observed that $\Theta_1(G)\leq d$
if and only if there are $d$ cliques in $G$ 
such that every edge of $G$ belongs to at least one of these cliques.
Kou et al. \cite{ksw} showed that deciding $\Theta_1(G)\leq k$ for a given graph is NP-complete.
Most of the research on $\Theta_p(G)$ focused on estimates for 
special graphs such as paths, trees, bounded degree graphs, complete bipartite graphs, see for instance \cite{cw,jkw,egr,mm}.

In the present paper we consider the classes 
$${\cal G}(d,p)=\{ G:\Theta_p(G)\leq d\}$$
of graphs for positive integers $d$ and $p$.
Clearly, ${\cal G}(d,p)$ is a hereditary class of graph, 
and can therefore be characterized by minimal forbidden induced subgraphs.
We give an upper bound on the order of minimal forbidden induced subgraphs for ${\cal G}(d,p)$.
In principle, for every choice of $d$ and $p$, 
this leads to a finite procedure that determines the complete list of minimal forbidden induced subgraphs for ${\cal G}(d,p)$.
Nevertheless, unless $d$ and $p$ are rather restricted, this procedure is impractical. 
For $p\in \{ d-1,d-2\}$, we provide more explicit results.

Considering the incidence vectors of the involved subsets of the ground set,
it is easy to see that some graph $G$
has a $p$-intersection representation $S:V(G)\to 2^U$ 
with $d=|U|$
if and only if 
there is a function $f:V(G)\to \{ 0,1\}^d$ 
such that distinct vertices $u$ and $v$ are adjacent in $G$ 
if and only if the dot product $f(u)\cdot f(v)$ of $f(u)$ and $f(v)$ is at least $p$.
We refer to such a function as a {\it binary dot product representation of dimension $d$ with threshold $p$}.
Clearly, $\Theta_p(G)$ is the minimum $d$ such that $G$ has a 
binary dot product representation of dimension $d$ with threshold $p$.
Dot product representations using real vectors and thresholds were studied for instance in \cite{fstz,klms}.

\section{Results}

Our first goal is an upper bound on the order of minimal forbidden induced subgraphs for ${\cal G}(d,p)$.
In fact, we consider slightly more general classes of graphs.

For a graph $G_0$ and a partition $V(G_0)=C\cup I$ of its vertex set, let ${\cal G}(G_0;C,I)$ 
denote the class of graphs that arise from $G$ by
\begin{itemize}
\item replacing
every vertex $u$ in $C$ by a possibly empty clique $C_u$, and
\item replacing every vertex $u$ in $I$ by a possibly empty independent set $I_u$.
\end{itemize}
Clearly, ${\cal G}(G_0;C,I)$ is a hereditary class of graphs.

If $d$ and $p$ are positive integers,
$G_{(d,p)}$ is the graph of order $2^d$ 
for which the bijection $f:V\left(G_{(d,p)}\right)\to \{ 0,1\}^d$ 
is a binary dot product representation of dimension $d$ with threshold $p$,
\begin{eqnarray*}
C_{(d,p)}&=&\left\{ {\bf x}\in \{ 0,1\}^d:{\bf x}\cdot {\bf x}\geq p\right\}\mbox{, and }\\
I_{(d,p)}&=& \{ 0,1\}^d\setminus C_{(d,p)}=\left\{ {\bf x}\in \{ 0,1\}^d:{\bf x}\cdot {\bf x}< p\right\},
\end{eqnarray*}
then a given graph $G$ has a binary $d$-dot representation with threshold $p$ if and only if $G$ belongs to
${\cal G}\left(G_{(d,p)};C_{(d,p)},I_{(d,p)}\right)$.
If for example $d=3$ and $p=2$, then $G_{(3,2)}$ is the disjoint union of a claw $K_{1,3}$ and four isolated vertices,
the set $C_{(3,2)}$ contains the vertices of the claw,
and the set $I_{(3,2)}$ contains the four isolated vertices, that is, all graphs that have a binary $3$-dot representation with threshold $2$ arise from $G_{(3,2)}$ by replacing the vertices of the claw by cliques, and the isolated vertices by independent sets.

We bound the order of minimal forbidden induced subgraphs for ${\cal G}(G_0;C,I)$.

\begin{theorem}\label{theorem1}
Let $G_0$ be a graph and let $V(G)=C\cup I$ be a partition of its vertex set. If $H$ is a minimal forbidden induced subgraph of ${\cal G}(G_0;C,I)$, then the order of $H$ is at most $4|C||I|+2|C|+2|I|+1$.
Specifically, ${\cal G}(G_0;C,I)=Forb({\cal F})$ for a finite set ${\cal F}$ of graphs.
\end{theorem}
{\it Proof:} 
First, we assume that there are at least $|I|+2$ vertices $u_1,\ldots,u_k$ of $H$ that are twins, 
that is, $N_H[u_1]=\ldots=N_H[u_k]$.
Since $H-u_k$ belongs to ${\cal G}(G_0;C,I)$, replacing the vertices $u$ in $C$ by suitable cliques $C_u$, 
and replacing the vertices $u$ in $I$ by suitable independent sets $I_u$,
results in $H-u_k$.
Since for every vertex $u$ in $I$, the set $I_u$ contains at most one of the vertices $u_1,\ldots,u_{k-1}$, 
and since $k-1\geq |I|+1$, there is some vertex $v$ in $C$ such that $u_i\in C_v$ for some $i\in [k]$. 
Replacing the vertices in $V(G)\setminus \{ v\}$ as before, and replacing the vertex $v$ by the clique $C_v\cup \{ u_k\}$ results in $H$, which is a contradiction. This implies that for every vertex $u$ of $H$, there are at most $|I|$ distinct further vertices of $H$ that have the same closed neighborhood as $u$.
Similarly, suitably exchanging the roles of $C$ and $I$ in the above argument, we obtain that for every vertex $u$ of $H$, there are at most $|C|$ distinct further vertices of $H$ that have the same (open) neighborhood as $u$.

Let $u^*$ be a vertex of $H$.
Let $H'=H-u^*$.
If there are at least $2|I|+3$ vertices of $H'$ that have the same closed neighborhood in $H'$, then at least $\lceil(2|I|+3)/2\rceil=|I|+2$ of these vertices have the same closed neighborhood in $H$, which is a contradiction.
Therefore, for every vertex $u$ of $H'$, there are at most $2|I|+1$ distinct further vertices of $H'$ that have the same closed neighborhood as $u$, and, similarly, for every vertex $u$ of $H'$, there are at most $2|C|+1$ distinct further vertices of $H'$ that have the same neighborhood as $u$.
Since $H'$ belongs to ${\cal G}(G_0;C,I)$, replacing the vertices $u$ in $C$ by suitable cliques $C'_u$, 
and replacing the vertices $u$ in $I$ by suitable independent sets $I'_u$, results in $H'$.
Since for every vertex $u\in C$, all vertices in $C'_u$ have the same closed neighborhood in $H'$, we have $|C'_u|\leq 2|I|+2$.
Since for every vertex $u\in I$, all vertices in $I'_u$ have the same neighborhood in $H'$, we have $|I'_u|\leq 2|C|+2$.
Altogether, we obtain
$n(H)=n(H')+1\leq |C|(2|I|+2)+|I|(2|C|+2)+1=
4|C||I|+2|C|+2|I|+1$. $\Box$

\bigskip

\noindent As a corollary, we obtain the desired upper bound 
on the order of minimal forbidden induced subgraphs for ${\cal G}(d,p)$.

\begin{corollary}\label{corollary1}
For positive integers $d$ and $p$, every minimal forbidden induced subgraph for ${\cal G}(d,p)$ has order at most $3\cdot 2^{d+1}+1$.
\end{corollary}
{\it Proof:} For the graph $G_{(d,p)}$, we have $|C|+|I|=|C_{(d,p)}|+|I_{(d,p)}|=2^d$, which implies $|C||I|\leq 2^d$, 
and hence $4|C||I|+2|C|+2|I|+1\leq 4\cdot 2^d+2\cdot 2^d+1$. 
$\Box$

\bigskip

\noindent The exponential dependence on $d$ in Corollary \ref{corollary1} is actually necessary.
This follows by choosing $p=\lceil d/2\rceil$ in the following result.

\begin{theorem}\label{theorem1b}
For positive integers $d$ and $p$, 
the graph $K_{1,{d\choose p}+1}$ is a minimal forbidden induced subgraph for ${\cal G}(d,p)$.
\end{theorem}
{\it Proof:} Let $k={d\choose p}+1$.
We need to show that $\Theta_p(K_{1,k})>d$,
and that $\Theta_p(G)\leq d$ for every proper induced subgraph $G$ of $K_{1,k}$.

For a contradiction, 
we assume that $f:V(K_{1,k})\to \{ 0,1\}^d$ 
is a binary dot product representation of dimension $d$ with threshold $p$
for $K_{1,k}$.
Let $I$ be the set of vertices of degree $1$ in $K_{1,k}$.
Since no vertex in $I$ is isolated, we have $f(u)\cdot f(u)\geq p$ for every $u\in I$.
Since no two vertices in $I$ are adjacent, this implies that $f(u)\not=f(v)$ for every two distinct vertices $u,v\in I$.
This implies the existence of $k$ distinct subsets $S_1,\ldots,S_k$ of $[d]$
with $|S_i|\geq p$ for $i\in [k]$,
and $|S_i\cap S_j|<p$ for distinct $i,j\in [k]$.
Possibly replacing sets $S_i$ of cardinality larger than $p$ 
with subsets $S_i'$ of cardinality exactly $p$,
this implies the existence of $k$ distinct subsets of order $p$ of $[d]$,
which is a contradiction, because $k>{d\choose p}$.
Therefore, $\Theta_p(K_{1,k})>d$.

Now let $G$ be a proper induced subgraph of $K_{1,k}$.
If all vertices of $G$ are isolated, 
then assigning to each vertex the all-$0$ vector of dimension $d$,
yields a binary dot product representation of dimension $d$ with threshold $p$ for $G$.
Hence, we may assume that $G$ is an induced subgraph of $K_{1,{d\choose p}}$.
Now, assigning to the universal vertex the all-$1$ vector of dimension $d$,
and to the remaining at most ${d\choose p}$ vertices distinct vectors from
$\left\{ {\bf x}\in \{ 0,1\}^d:{\bf x}\cdot {\bf x}=p\right\}$,
yields a binary dot product representation of dimension $d$ with threshold $p$ for $G$,
which completes the proof. $\Box$

\bigskip

\noindent We proceed to a more explicit result for $p=d-1$.
In order to reduce the number of different minimal forbidden induced subgraphs further,
we exclude isolated vertices. 

\begin{theorem}\label{theorem2}
Let $G$ be a graph without isolated vertices, 
and let $d$ be a positive integer at least $2$. 
The following statements are equivalent.
\begin{enumerate}[(i)]
\item $\Theta_{d-1}(G)\leq d$.
\item $G$ arises by replacing the vertices of a claw $K_{1,d}$ by possibly empty cliques.
\item $G$ belongs to $Forb\left(\left\{ \bar{K}_{d+1},P_4,K_1\cup P_3,C_4\right\}\right)$.
\end{enumerate}
\end{theorem}
{\it Proof:} The observations before Theorem \ref{theorem1} imply the equivalence of (i) and (ii).
Since it is easy to see that (ii) implies (iii), it remains to show that (iii) implies (ii).
Therefore, let $G$ be a graph without isolated vertices that belongs to $Forb(\{ \bar{K}_{d+1},P_4,K_1\cup P_3,C_4\})$.
If $G$ is not connected, then, 
since $G$ is $(K_1\cup P_3)$-free, 
all components of $G$ are cliques, 
and, since $G$ is $\bar{K}_{d+1}$-free, $G$ has at most $d$ components.
Altogether, in this case it follows that $G$ is as in (ii).
Now let $G$ be connected.
Clearly, we may assume that $G$ is not a clique.
Since $G$ is $P_4$-free, 
there is a non-trivial partition $V(G)=V_1\cup V_2$ such that $G$ contains all edges between $V_1$ and $V_2$.
Since $G$ is $C_4$-free, we may assume that $V_1$ is a clique.
Furthermore, we may assume that the partition $V(G)=V_1\cup V_2$ is chosen such that 
$|V_1|$ is as large as possible. 
Since $G$ is not a clique, also $G[V_2]$ is not a clique.
If $G[V_2]$ is connected, then using $\{ P_4,C_4\}$-freeness as above implies the existence of a universal vertex $u$ of $G[V_2]$, and the partition $V(G)=(V_1\cup \{ u\})\cup (V_2\setminus \{ u\})$ 
contradicts the choice of the partition $V(G)=V_1\cup V_2$.
Hence, $G[V_2]$ is not connected.
Since $G$ is $(K_1\cup P_3)$-free, all components of $G[V_2]$ are cliques, and, 
since $G$ is $\bar{K}_{d+1}$-free, $G[V_2]$ has at most $d$ components.
Altogether, also in this case it follows that $G$ is as in (ii),
which completes the proof. $\Box$

\bigskip

\noindent We proceed to the case $p=d-2$.
In order to reduce the number of different minimal forbidden induced subgraphs further,
we exclude isolated as well as universal vertices this time.

For a positive integer $d$ at least $3$, 
let $G^*_{(d,d-2)}$ be the graph of order $d+{d\choose 2}$ 
for which the bijection 
$$f:V\left(G^*_{(d,d-2)}\right)\to \left\{ {\bf x}\in\{ 0,1\}^d:d-2\leq {\bf x}\cdot {\bf x}\leq d-1\right\}$$ 
is a binary dot product representation of dimension $d$ with threshold $d-2$.
Note that $G^*_{(d,d-2)}$ is a split graph
whose vertex set is partitioned into a clique $C$ of order $d$
and an independent set $I$ of order ${d\choose 2}$
such that every vertex in $I$ has degree $2$,
and for every pair $c_1$ and $c_2$ of distinct vertices in $C$,
there is a vertex $i$ in $I$ that is adjacent to $c_1$ and $c_2$.

For a graph $G$, let the {\it reduction ${\cal R}(G)$} of $G$ arise from $G$ by 
identifying all pairs of vertices that are twins.

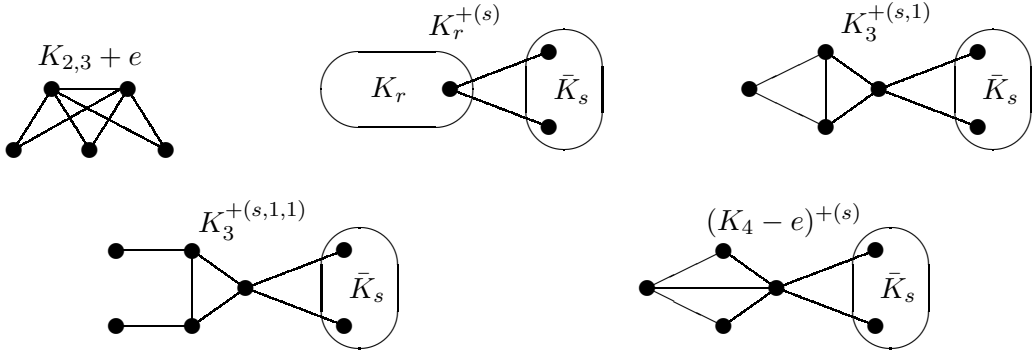
\begin{figure}[H]
\begin{center}
$\mbox{}$\hfill
\unitlength 1mm 
\linethickness{0.4pt}
\ifx\plotpoint\undefined\newsavebox{\plotpoint}\fi 
\begin{picture}(21,11)(0,0)
\put(0,0){\circle*{2}}
\put(10,0){\circle*{2}}
\put(20,0){\circle*{2}}
\put(5,8){\circle*{2}}
\put(15,8){\circle*{2}}
\multiput(0,0)(.033557047,.053691275){149}{\line(0,1){.053691275}}
\multiput(5,8)(.033557047,-.053691275){149}{\line(0,-1){.053691275}}
\multiput(10,0)(.033557047,.053691275){149}{\line(0,1){.053691275}}
\multiput(15,8)(.033557047,-.053691275){149}{\line(0,-1){.053691275}}
\multiput(20,0)(-.06302521,.033613445){238}{\line(-1,0){.06302521}}
\put(5,8){\line(1,0){9}}
\multiput(14,8)(-.058823529,-.033613445){238}{\line(-1,0){.058823529}}
\put(10,12){\makebox(0,0)[cc]{$K_{2,3}+e$}}
\end{picture}
\hfill
\unitlength 1mm 
\linethickness{0.4pt}
\ifx\plotpoint\undefined\newsavebox{\plotpoint}\fi 
\begin{picture}(37,16)(0,0)
\put(10,8){\oval(20,10)[]}
\put(17,8){\circle*{2}}
\put(9,8){\makebox(0,0)[cc]{$K_r$}}
\put(30,3){\circle*{2}}
\put(30,13){\circle*{2}}
\multiput(30,13)(-.087248322,-.033557047){149}{\line(-1,0){.087248322}}
\multiput(17,8)(.087248322,-.033557047){149}{\line(1,0){.087248322}}
\put(33,8){\makebox(0,0)[cc]{$\bar{K}_s$}}
\put(32,8){\oval(10,16)[]}
\put(19,17){\makebox(0,0)[cc]{$K_r^{+(s)}$}}
\end{picture}
\hfill
\unitlength 1mm 
\linethickness{0.4pt}
\ifx\plotpoint\undefined\newsavebox{\plotpoint}\fi 
\begin{picture}(37,16)(0,0)
\put(17,8){\circle*{2}}
\put(30,3){\circle*{2}}
\put(30,13){\circle*{2}}
\multiput(30,13)(-.087248322,-.033557047){149}{\line(-1,0){.087248322}}
\multiput(17,8)(.087248322,-.033557047){149}{\line(1,0){.087248322}}
\put(33,8){\makebox(0,0)[cc]{$\bar{K}_s$}}
\put(32,8){\oval(10,16)[]}
\put(10,3){\circle*{2}}
\put(10,13){\circle*{2}}
\multiput(17,8)(-.046979866,.033557047){149}{\line(-1,0){.046979866}}
\put(10,13){\line(0,-1){10}}
\multiput(10,3)(.046979866,.033557047){149}{\line(1,0){.046979866}}
\put(18,17){\makebox(0,0)[cc]{$K_3^{+(s,1)}$}}
\put(0,8){\circle*{2}}
\put(10,13){\line(-2,-1){10}}
\put(0,8){\line(2,-1){10}}
\end{picture}
\hfill $\mbox{}$\\[10mm]
$\mbox{}$\hfill
\unitlength 1mm 
\linethickness{0.4pt}
\ifx\plotpoint\undefined\newsavebox{\plotpoint}\fi 
\begin{picture}(37,16)(0,0)
\put(17,8){\circle*{2}}
\put(30,3){\circle*{2}}
\put(30,13){\circle*{2}}
\multiput(30,13)(-.087248322,-.033557047){149}{\line(-1,0){.087248322}}
\multiput(17,8)(.087248322,-.033557047){149}{\line(1,0){.087248322}}
\put(33,8){\makebox(0,0)[cc]{$\bar{K}_s$}}
\put(32,8){\oval(10,16)[]}
\put(10,3){\circle*{2}}
\put(10,13){\circle*{2}}
\multiput(17,8)(-.046979866,.033557047){149}{\line(-1,0){.046979866}}
\put(10,13){\line(0,-1){10}}
\multiput(10,3)(.046979866,.033557047){149}{\line(1,0){.046979866}}
\put(0,3){\circle*{2}}
\put(0,3){\line(1,0){9}}
\put(18,17){\makebox(0,0)[cc]{$K_3^{+(s,1,1)}$}}
\put(0,13){\circle*{2}}
\put(0,13){\line(1,0){10}}
\end{picture}
\hfill
\unitlength 1mm 
\linethickness{0.4pt}
\ifx\plotpoint\undefined\newsavebox{\plotpoint}\fi 
\begin{picture}(37,16)(0,0)
\put(17,8){\circle*{2}}
\put(30,3){\circle*{2}}
\put(30,13){\circle*{2}}
\multiput(30,13)(-.087248322,-.033557047){149}{\line(-1,0){.087248322}}
\multiput(17,8)(.087248322,-.033557047){149}{\line(1,0){.087248322}}
\put(33,8){\makebox(0,0)[cc]{$\bar{K}_s$}}
\put(32,8){\oval(10,16)[]}
\put(10,3){\circle*{2}}
\put(10,13){\circle*{2}}
\multiput(17,8)(-.046979866,.033557047){149}{\line(-1,0){.046979866}}
\multiput(10,3)(.046979866,.033557047){149}{\line(1,0){.046979866}}
\put(18,17){\makebox(0,0)[cc]{$(K_4-e)^{+(s)}$}}
\put(0,8){\circle*{2}}
\put(10,13){\line(-2,-1){10}}
\put(0,8){\line(2,-1){10}}
\put(0,8){\line(1,0){17}}
\end{picture}
\hfill$\mbox{}$
\caption{
The graphs
$K_{2,3}+e$,
$K_r^{+(s)}$,
$K_3^{+(s,1)}$,
$K_3^{+(s,1,1)}$, and
$(K_4-e)^{+(s)}$.}\label{fig1}
\end{center}
\end{figure}

Let $K_{2,3}+e$ arise by adding an edge between the two vertices of degree $3$ of $K_{2,3}$.
Let $K_n-e$ arise by removing one edge from $K_n$.
Let $K_r^{+(s)}$ arise by adding $s$ pendant vertices to one vertex of $K_r$.
Note that $K_2^{+(t-1)}=K_{1,t}$.
Let $K_3^{+(s,1)}$ arise from $K_3$
by adding $s$ pendant vertices to one vertex of $K_3$,
and adding one more vertex that is adjacent to the other two vertices of $K_3$.
Let $K_3^{+(s,1,1)}$ arise from $K_3$
by adding $s$ pendant vertices to one vertex of $K_3$,
and adding one pendant vertex to each of the other two vertices of $K_3$.
Let $(K_4-e)^{+(s)}$ arise by adding $s$ pendant vertices to a vertex of degree $3$ in $K_4-e$.
See Figure \ref{fig1} for illustrations of these graphs.

\begin{theorem}\label{theorem3}
Let $G$ be a graph without isolated or universal vertices,
and let $d$ be a positive integer at least $3$. 
The following statements are equivalent.
\begin{enumerate}[(i)]
\item $\Theta_{d-2}(G)\leq d$.
\item $G$ arises by replacing the vertices of $G^*_{(d,d-2)}$ by possibly empty cliques.
\item ${\cal R}(G)$ belongs to $Forb({\cal F})$ for 
\begin{eqnarray*}
{\cal F} & = & 
\Big\{ 
2K_2,
C_4,
C_5,
K_{2,3}+e,
K_5-e
\Big\}\\
&& 
\cup \Big\{ 
K_3^{+(d-2,1)},
K_3^{+(d-2,1,1)},
(K_4-e)^{+(d-2)}
\Big\}\\
&& 
\cup \Big\{ 
K_i^{+(d+1-i)}:i=2\mbox{ or }4\leq i\leq d+1
\Big\}\\
&& 
\cup \Big\{ 
K_2\cup \bar{K}_{{d-1\choose 2}+1},
K_3\cup \bar{K}_{{d-2\choose 2}+1},
P_4\cup \bar{K}_{{d-2\choose 2}+1}
\Big\}\\
&&
\cup \Big\{ 
K_3^{+(1,1)}\cup \bar{K}_{{d-3\choose 2}+1},
K_3^{+(1,1,1)}\cup \bar{K}_{{d-3\choose 2}+1},
(K_4-e)\cup \bar{K}_{{d-3\choose 2}+1}
\Big\}\\
&& 
\cup \Big\{ 
K_i\cup \bar{K}_{{d-i\choose 2}+1}:i=0\mbox{ or }4\leq i\leq d
\Big\}.
\end{eqnarray*}
\end{enumerate}
\end{theorem}
{\it Proof:} Let $G$ satisfy (i). 
Let $f$ be a binary dot representation of dimension $d$ and threshold $d-2$ for $G$.
Since $G$ has no isolated vertex, for every vertex $u$ of $G$, 
the vector $f(u)$ contains at most two $0$-entries. 
Since $G$ has no universal vertex, for every vertex $u$ of $G$, 
the vector $f(u)$ contains at least one $0$-entry.
Altogether, for every vertex $u$ of $G$,
the vector $f(u)$ belongs to $\left\{ {\bf x}\in\{ 0,1\}^d:d-2\leq {\bf x}\cdot {\bf x}\leq d-1\right\}$.

Since in $\{ 0,1\}^d$, there are disjoint sets $C$ and $I$ such that 
\begin{itemize}
\item $C$ contains $d$ vectors with exactly one $0$-entry, every two of which have dot product $d-2$,
\item $I$ contains ${d\choose 2}$ vectors with exactly two $0$-entries, every two of which have dot product at most $d-3$, and
\item for every two vectors ${\bf c}_1$ and ${\bf c}_2$ in $C$,
there is a vector ${\bf i}$ in $I$ such that ${\bf c}_1\cdot {\bf i},{\bf c}_2\cdot {\bf i}\geq d-2$,
\end{itemize}
it follows that (i) implies (ii).
Conversely, the definition of $G^*_{(d,d-2)}$ easily implies that (ii) implies (i),
that is, (i) and (ii) are equivalent.

\bigskip

\noindent It is tedious yet not difficult to show that (i) implies (iii).
We give details for the graphs in
$$\left\{ 2K_2, C_5,K_{2,3}+e,K_3^{+(1,1,1)}\cup \bar{K}_{{d-3\choose 2}+1}\right\} 
\cup \left\{ K_i^{+(d+1-i)}:4\leq i\leq d\right\},$$
and leave the details for the remaining graphs to the reader.
Therefore, let $G$ be as in (i).
Let $f:V(G)\to \{ 0,1\}^d$ be a binary dot representation of dimension $d$ and threshold $d-2$ for $G$.
As observed above, 
$f(u)\in \left\{ {\bf x}\in\{ 0,1\}^d:d-2\leq {\bf x}\cdot {\bf x}\leq d-1\right\}$
for every vertex $u$ of $G$.
Note that ${\cal R}(G)$ is an induced subgraph of $G$.
Let 
\begin{eqnarray*}
C &=& \left\{u\in V\left({\cal R}(G)\right):f(u)\cdot f(u)=d-1\right\}\mbox{ and}\\
I &=& V\left({\cal R}(G)\right)\setminus C= \left\{u\in V\left({\cal R}(G)\right):f(u)\cdot f(u)=d-2\right\}.
\end{eqnarray*}
Note that 
$C$ is a clique of $G$ of order at most $d$, 
$I$ is an independent set of $G$ of order at most ${d\choose 2}$,
every vertex in $I$ is adjacent to at most two vertices in $C$, and
for every two vertices in $C$, there is at most one vertex in $I$ that is adjacent to both.
In particular, all vertices in $I$ are simplicial.

If ${\cal R}(G)$ contains $2K_2$ as an induced subgraph,
then one of the components of $2K_2$ does not intersection $C$, 
because $C$ is a clique.
Now, this implies the existence of an edge within $I$,
which is a contradiction.
Hence, ${\cal R}(G)$ is $2K_2$-free.
If ${\cal R}(G)$ contains $C_5$ as an induced subgraph,
then all five non-simplicial vertices of $C_5$ belong to $C$,
which is a contradiction.
Hence, ${\cal R}(G)$ is $C_5$-free.
If ${\cal R}(G)$ contains $K_{2,3}+e$ as an induced subgraph,
then the two non-simplicial vertices of $K_{2,3}+e$ belong to $C$.
Since $C$ is a clique, at least two of the three vertices of degree $2$ of $K_{2,3}+e$ belong to $I$.
This implies the existence of two vertices in $I$ having the same two neighbors in $C$,
which is a contradiction.
Hence, ${\cal R}(G)$ is $(K_{2,3}+e)$-free.

If ${\cal R}(G)$ contains $K_3^{+(1,1,1)}\cup \bar{K}_{{d-3\choose 2}+1}$ as an induced subgraph,
then the three non-simplicial vertices of degree $3$, say $x_1$, $x_2$, and $x_3$, 
of $K_3^{+(1,1,1)}\cup \bar{K}_{{d-3\choose 2}+1}$ belong to $C$,
and the ${d-3\choose 2}+1$ isolated vertices of $K_3^{+(1,1,1)}\cup \bar{K}_{{d-3\choose 2}+1}$ belong to $I$.
Since no two vertices in $C \cup I$ are twins in $G$, the function $f$ is injective on $C \cup I$. 
This implies the existence of 
${d-3\choose 2}+1$
distinct vectors in 
$$U'=\left\{ {\bf x}\in\{ 0,1\}^d:{\bf x}\cdot {\bf x}=d-2\right\}
\setminus \bigcup_{j=1}^3\left\{ {\bf x}\in\{ 0,1\}^d:{\bf x}\cdot f(x_j)=d-2\right\},$$ 
where $f(x_1)$, $f(x_2)$, and $f(x_3)$ are three distinct vectors in 
$\left\{ {\bf x}\in\{ 0,1\}^d:{\bf x}\cdot {\bf x}=d-1\right\}$.
Since $U'$ contains exactly ${d\choose 2}-(d-1)-(d-2)-(d-3)={d-3\choose 2}$ elements, 
we obtain a contradiction.
Hence, ${\cal R}(G)$ is $\left( K_3^{+(1,1,1)}\cup \bar{K}_{{d-3\choose 2}+1}\right)$-free.

If ${\cal R}(G)$ contains $K_i^{+(d+1-i)}$ as an induced subgraph for some $4\leq i\leq d$,
then the $i$ vertices of degree at least $3$ of $K_i^{+(d+1-i)}$ belong to $C$,
and the $d+i-1$ vertices of degree $1$ of $K_i^{+(d+1-i)}$ belong to $I$.
Let $x_1$ be the vertex of degree at least $i$ of $K_i^{+(d+1-i)}$,
and let $x_2,\ldots,x_i$ be the vertices of degree $i-1$ of $K_i^{+(d+1-i)}$.
Since $f$ is injective on $C \cup I$, 
we obtain the existence of $d+1-i$ distinct vectors in 
$$U''=\left\{ {\bf x}\in\{ 0,1\}^d:{\bf x}\cdot {\bf x}=d-2\mbox{ and }{\bf x}\cdot f(x_1)=d-2\right\}
\setminus \bigcup_{j=2}^i\left\{ {\bf x}\in\{ 0,1\}^d:{\bf x}\cdot f(x_j)=d-2\right\},$$
where $f(x_1),\ldots,f(x_i)$ are $i$ distinct vectors in 
$\left\{ {\bf x}\in\{ 0,1\}^d:{\bf x}\cdot {\bf x}=d-1\right\}$. 
Since $U''$ contains exactly $d-i$ elements,
we obtain a contradiction.
Hence, ${\cal R}(G)$ is $\left( K_i^{+(d+1-i)}\right)$-free.

\bigskip

\noindent We proceed to show that (iii) implies (ii).
Therefore, let $G$ be a graph without isolated or universal vertices 
such that (iii) holds.
Let $H={\cal R}(G)$.
Let $X$ be a maximum clique of $H$.
If possible, we choose $X$ in such a way that $d_H(u)\geq |X|$ for every $u\in X$.
Let $k=|X|$.
Let $Y$ denote the set of vertices in $V(H)\setminus X$ that are not isolated in $H$.
Let $Z=V(H)\setminus (X\cup Y)$.
We consider different cases.

\bigskip

\noindent {\bf Case 1} $k\geq 4$.

\bigskip

\noindent Since $H$ is $(K_5-e)$-free, every vertex in $Y$ has at most two neighbors in $X$.
For a contradiction, we assume that $Y\cup Z$ is not independent.
This implies that there is some edge $y_1y_2$ with $y_1,y_2\in Y$.
Since $H$ is $2K_2$-free, each but at most one vertex from $X$ has a neighbor in $\{ y_1,y_2\}$.
Since $y_1$ and $y_2$ both have at most two neighbors in $X$,
this implies that there are two vertices $x_1,x_2\in X$ with 
$x_1y_1,x_2y_2\in E(H)$ and $x_1y_2,x_2y_1\not\in E(H)$.
Now, $x_1y_1y_2x_2x_1$ is a $C_4$,
which is a contradiction.
Hence, $Y\cup Z$ is independent.
Since $H$ is $(K_{2,3}+e)$-free, for every two vertices $x_1$ and $x_2$ in $X$,
there is at most one vertex $y$ in $Y$ that is adjacent to $x_1$ and $x_2$.
Since $H$ is $K_k^{+(d-k+1)}$-free, 
every vertex in $H$ is adjacent to at most $d-k$ vertices of degree $1$.
Since $H$ is $\left(K_k\cup \bar{K}_{{d-k\choose 2}+1}\right)$-free,
$H$ has at most ${d-k\choose 2}$ isolated vertices.

Let $X=\{ x_1,\ldots,x_k\}$.
We will now describe a binary dot representation of dimension $d$ with threshold $d-2$ for $H$,
which will imply that $H$, and hence also $G$, satisfies (ii):
\begin{itemize}
\item For $i\in [k]$, assign to $x_i$ the vector ${\bf 1}-{\bf e}_i$, 
where ${\bf 1}$ is the all-$1$ vector of dimension $d$,
and ${\bf e}_i$ is the $i$-th unit vector of dimension $d$.
\item For every vertex $v$ in $Y$ that is adjacent to $x_i$ and $x_j$ for distinct $i$ and $j$ in $[k]$, 
assign to $v$ the vector ${\bf 1}-{\bf e}_i-{\bf e}_j$.
\item For $i\in [k]$, assign to the at most $d-k$ neighbors of $x_i$ of degree $1$ 
distinct vectors from $\{ {\bf 1}-{\bf e}_i-{\bf e}_j:k+1\leq j\leq d\}$.
\item Assign to the at most ${d-k\choose 2}$ isolated vertices
distinct vectors from $\{ {\bf 1}-{\bf e}_i-{\bf e}_j:k+1\leq i<j\leq d\}$.
\end{itemize}

\bigskip

\noindent Let $X=\{ x_1,\ldots,x_k\}$.

\bigskip

\noindent {\bf Case 2} $k=3$ and $d_H(u)\geq 3$ for every $u\in X$.

\bigskip

\noindent By the choice of $X$, every vertex in $Y$ has at most two neighbors in $X$.

For a contradiction, we assume that $Y\cup Z$ is not independent.
This implies that there is some edge $y_1y_2$ with $y_1,y_2\in Y$.
Since $H$ is $\{ 2K_2,C_4\}$-free, 
we may assume, by symmetry, that $N_H(y_1)\cap X=\{ x_1,x_2\}$.
If $y_2$ is adjacent to $x_3$, 
then either $y_1y_2x_3x_1y_1$ or $y_1y_2x_3x_2y_1$ is a $C_4$,
which is a contradiction.
Hence, $y_2$ is not adjacent to $x_3$. 
Since $x_3$ has degree at least $3$, 
there is a neighbor $y_3$ of $x_3$ that is distinct from $x_1$ and $x_2$.
Since $H$ is $2K_2$-free, $y_3$ has a neighbor in $\{ y_1,y_2\}$.
Since $y_3$ has at most two neighbors in $X$,
we may assume that $x_1$ is not adjacent to $y_3$.
Now 
either $y_1y_3x_3x_1y_1$,
or $y_2y_3x_3x_1y_2$,
or $y_2y_3x_3x_1y_1y_2$
is $C_4$ or $C_5$, which is a contradiction.
Hence, $Y\cup Z$ is independent.

Since $H$ is $(K_{2,3}+e)$-free, for every two vertices $x_1$ and $x_2$ in $X$,
there is at most one vertex $y$ in $Y$ that is adjacent to $x_1$ and $x_2$.

For a contradiction, we assume that $x_1$ is adjacent to $d-2$ vertices of degree $1$.
Since $x_2$ and $x_3$ have degree at least $3$, there are 
not necessarily distinct neighbors $y_2$ of $x_2$ and $y_3$ of $x_3$
that do not belong to $X$.
If $y_2=y_3$, then $H$ contains $K_3^{+(d-2,1)}$, which is a contradiction.
Hence, $y_2\not=y_3$.
We may assume that $x_2$ is not adjacent to $y_3$, and $x_3$ is not adjacent to $y_2$.
Since $H$ is $K_{1,d}$-free, $x_1$ is not adjacent to $y_2$ or $y_3$.
Now, $H$ contains $K_3^{+(d-2,1,1)}$, which is a contradiction.
Hence, by symmetry,
each vertex in $X$ is adjacent to at most $d-3$ vertices of degree $1$.

Since $H$ is $\left\{ K_3^{+(1,1)}\cup \bar{K}_{{d-3\choose 2}+1},
K_3^{+(1,1,1)}\cup \bar{K}_{{d-3\choose 2}+1}\right\}$-free,
$H$ has at most ${d-3\choose 2}$ isolated vertices.
At this point, we can complete the proof as in Case 1 setting $k=3$.

\bigskip

\noindent {\bf Case 3} $k=3$ and every triangle in $H$ contains a vertex of degree $2$.

\bigskip

\noindent If $y_1y_2$ is an edge between two vertices of $Y$,
then we may assume as above that $N_H(y_1)\cap X=\{ x_1,x_2\}$.
Now $y_1x_1x_2y_1$ is a triangle in which every vertex has degree at least $3$,
which is a contradiction. Hence, $Y\cup Z$ is independent.

First, we assume that $x_3$ is the only common neighbor of $x_1$ and $x_2$.
Since $H$ is $K_2^{+(d-1)}$-free, 
each vertex in $\{ x_1,x_2\}$ is adjacent to at most $d-2$ vertices of degree $1$.
Since $H$ is $\left( K_3\cup \bar{K}_{{d-2\choose 2}+1}\right)$-free,
$H$ has at most ${d-2\choose 2}$ isolated vertices.
At this point, we can complete the proof as in Case 1 setting $k=2$.

Next, we assume that $x_1$ and $x_2$ have a second common neighbor $y_3$ distinct from $x_3$.
Since $H$ is $(K_{2,3}+e)$-free, $x_1$ and $x_2$ have exactly two common neighbors.
Since $H$ is $(K_4-e)^{+(d-2)}$-free, 
each vertex in $\{ x_1,x_2\}$ is adjacent to at most $d-3$ vertices of degree $1$.
Since $H$ is $\left((K_4-e)\cup \bar{K}_{{d-3\choose 2}+1}\right)$-free,
$H$ has at most ${d-3\choose 2}$ isolated vertices.
At this point, we can complete the proof as in Case 1 setting $k=3$.

\bigskip

\noindent {\bf Case 4} $k=2$ and $d_H(u)\geq 2$ for every $u\in X$.

\bigskip

\noindent  By the choice of $X$, every vertex in $Y$ has at most one neighbor in $X$.
Since $H$ is $\{ 2K_2,C_4\}$-free, this implies that $Y\cup Z$ is independent.
Since $H$ is $K_2^{+(d-1)}$-free, 
each vertex in $\{ x_1,x_2\}$ is adjacent to at most $d-2$ vertices of degree $1$.
Since $H$ is $\left( P_4\cup \bar{K}_{{d-2\choose 2}+1}\right)$-free,
$H$ has at most ${d-2\choose 2}$ isolated vertices.
At this point, we can complete the proof as in Case 1 setting $k=2$.

\bigskip

\noindent {\bf Case 5} $H$ has at most one vertex of degree at least $2$.

\bigskip

\noindent  Since $H$ is $\left\{K_2^{+(d-1)},K_2\cup \bar{K}_{{d-1\choose 2}+1}\right\}$-free,
$H$ is a subgraph of $K_{1,d-1}\cup \bar{K}_{{d-1\choose 2}}$.
At this point, we can complete the proof as in Case 1 setting $k=1$. $\Box$

\bigskip

\noindent  {\bf Acknowledgment} 
J.W.C. Pinto and J.L. Szwarcfiter were partially supported by CAPES and CNPq.


\begin{thebibliography}{}

\bibitem{cw}
M.S. Chung and D.B. West, The $p$-intersection number of a complete bipartite graph and orthogonal double coverings of a clique, Combinatorica 14 (1994) 453-461.
\bibitem{egr}
N. Eaton, R. Gould, and V. R\"{o}dl, On p-intersection representations, J. Graph Theory 21 (1996) 377-392.
\bibitem{egp}
P. Erd\H{o}s, A. Goodman, and L. P\'{o}sa, The representation of a graph by set intersections, {Can. J. Math.} 18 (1966) 106-112.
\bibitem{fstz}
C.M. Fiduccia, E.R. Scheinermann, A. Trenk, and J.S. Zito, Dot product representations of graphs, {Discrete Math.} 181 (1998) 113-138.
\bibitem{jkw}
M.S. Jacobson, A. K\'{e}zdy, and D.B. West, The $2$-intersection number of paths and bounded-degree trees, J. Graph Theory 4 (1995) 461-469.
\bibitem{jms}
M.S. Jacobson, F.R. McMorris, and E.R. Scheinerman, General results on tolerance intersection graphs, J. Graph Theory 15 (1991) 573-577.
\bibitem{klms}
R.J. Kang, L. Lov\'{a}sz, T. M\"{u}ller, and E.R. Scheinerman, Dot product representations of planar graphs, Electronic J. Combin. 18 (2011) $\#$ P216.
\bibitem{ksw}
L.T. Kou, L.J. Stokmeyer, and C.K. Wong, Covering edges by cliques with regard to keyword conflicts and intersection graphs, Comm. ACM 21 (1978) 135-139.
\bibitem{mm}
T.A. McKee and F.R. McMorris, Topics in intersection graph theory, SIAM monographs on discrete mathematics and applications, 1999.
\end{thebibliography}
\end{document}